\magnification=\magstep1

\def\to{\ \longrightarrow\ }

\def\nl{\hfill\break}

\def\hexnumber#1{\ifcase#1 0\or 1\or 2\or 3\or 4\or 5\or 6\or 7\or 8\or
 9\or A\or B\or C\or D\or E\or F\fi}
%
%
\font\twelvemsa=msam10 scaled 1200   
\font\tenmsa=msam10                  
\font\ninemsa=msam9            \font\sevenmsa=msam7
\font\sixmsa=msam6             \font\fivemsa=msam5
%
%
\newfam\msafam                 \textfont\msafam=\tenmsa
\scriptfont\msafam=\sevenmsa   \scriptscriptfont\msafam=\fivemsa
\edef\hexa{\hexnumber\msafam}        
\def\msa{\fam\msafam\tenmsa}         
%
%
\font\twelvemsb=msbm10 scaled 1200   
\font\tenmsb=msbm10                  
\font\ninemsb=msbm9            \font\sevenmsb=msbm7
\font\sixmsb=msbm6             \font\fivemsb=msbm5
%
\newfam\msbfam                 \textfont\msbfam=\tenmsb       
\scriptfont\msbfam=\sevenmsb   \scriptscriptfont\msbfam=\fivemsb
\edef\hexb{\hexnumber\msbfam}        
\def\msb{\fam\msbfam\tenmsb}         
%
%
\font\twelveeufm=eufm10 scaled 1200  
\font\teneufm=eufm10                 
\font\nineeufm=eufm9           \font\seveneufm=eufm7
\font\sixeufm=eufm6            \font\fiveeufm=eufm5
%
\newfam\eufmfam                \textfont\eufmfam=\teneufm
\scriptfont\eufmfam=\seveneufm \scriptscriptfont\eufmfam=\fiveeufm
\edef\hexf{\hexnumber\eufmfam}      
\def\frak{\fam\eufmfam\teneufm}     
%
%
%
\font\twelverm=cmr10 scaled 1200    
\font\ninerm=cmr9                   
\font\sixrm=cmr6   
%
\font\twelvei=cmmi10 scaled 1200    
\font\ninei=cmmi9                   
\font\sixi=cmmi6  
%
\font\twelvesy=cmsy10 scaled 1200   
\font\ninesy=cmsy9                  
\font\sixsy=cmsy6  
%
\font\twelvebf=cmbx10 scaled 1200   
\font\ninebf=cmbx9                  
\font\sixbf=cmbx6  
%
%
\font\twelveit=cmti10 scaled 1200   
\font\nineit=cmti9                  
%
\font\twelvesl=cmsl10 scaled 1200   
\font\ninesl=cmsl9                  
%
\font\twelvett=cmtt10 scaled 1200   
\font\ninett=cmtt9                  
%
%
%
%
\def\small{%
%
%
\textfont0=\ninerm \scriptfont0=\sixrm \scriptscriptfont0=\fiverm
\def\rm{\fam0\ninerm}        
%
%
\textfont1=\ninei \scriptfont1=\sixi \scriptscriptfont1=\fivei
%
%
\textfont2=\ninesy \scriptfont2=\sixsy \scriptscriptfont2=\fivesy
%
%
\textfont3=\tenex \scriptfont3=\tenex \scriptscriptfont3=\tenex
%
%
\textfont\bffam=\ninebf \scriptfont\bffam=\sixbf
\scriptscriptfont\bffam=\fivebf \def\bf{\fam\bffam\ninebf}%
%
%
\textfont\itfam=\nineit \def\it{\fam\itfam\nineit}%
\textfont\slfam=\ninesl \def\sl{\fam\slfam\ninesl}%
\textfont\ttfam=\ninett \def\tt{\fam\ttfam\ninett}%
%
%
%
\textfont\msafam=\ninemsa \scriptfont\msafam=\sixmsa
\scriptscriptfont\msafam=\fivemsa \def\msa{\fam\msafam\ninemsa}%
%
%
\textfont\msbfam=\ninemsb \scriptfont\msbfam=\sixmsb
\scriptscriptfont\msbfam=\fivemsb \def\msb{\fam\msbfam\ninemsb}%
%
%
\textfont\eufmfam=\nineeufm  \scriptfont\eufmfam=\sixeufm
\scriptscriptfont\eufmfam=\fiveeufm \def\frak{\fam\eufmfam\nineeufm}%
%
%
%
\normalbaselineskip=11pt
\setbox\strutbox=\hbox{\vrule height8pt depth3pt width0pt}%
%
%
\normalbaselines\rm}    
%
%
%
%
\def\large{%
\textfont0=\twelverm \scriptfont0=\ninerm \scriptscriptfont0=\sevenrm
\def\rm{\fam0\twelverm}%
\textfont1=\twelvei \scriptfont1=\ninei \scriptscriptfont1=\seveni
\textfont2=\twelvesy \scriptfont2=\ninesy \scriptscriptfont2=\sevensy
\textfont3=\tenex \scriptfont3=\tenex \scriptscriptfont3=\tenex
\textfont\bffam=\twelvebf \scriptfont\bffam=\ninebf
\scriptscriptfont\bffam=\sevenbf \def\bf{\fam\bffam\twelvebf}%
\textfont\itfam=\twelveit \def\it{\fam\itfam\twelveit}%
\textfont\slfam=\twelvesl \def\sl{\fam\slfam\twelvesl}%
\textfont\ttfam=\twelvett \def\tt{\fam\ttfam\twelvett}%
\textfont\msafam=\twelvemsa \scriptfont\msafam=\ninemsa
\scriptscriptfont\msafam=\sevenmsa \def\msa{\fam\msafam\twelvemsa}         
\textfont\msbfam=\twelvemsb \scriptfont\msbfam=\ninemsb
\scriptscriptfont\msbfam=\sevenmsb \def\msb{\fam\msbfam\twelvemsb}         
\textfont\eufmfam=\twelveeufm  \scriptfont\eufmfam=\nineeufm
\scriptscriptfont\eufmfam=\seveneufm \def\frak{\fam\eufmfam\teneufm}
\normalbaselineskip=15pt
\setbox\strutbox=\hbox{\vrule height11pt depth4pt width0pt}%
\normalbaselines\rm}%
%

%

%
\mathchardef\plussquare="0\hexa01
\mathchardef\nge="3\hexb0B
\mathchardef\maltesecross="0\hexa7A
\mathchardef\del="0\hexf01
%

%

\input epsf
\input colordvi

\overfullrule=0pt

\font\npt=cmr9

\font\secfont=cmbx10

\font\nam=cmr8
\font\aff=cmti8

\mathchardef\square="0\hexa03
\def\qed{\hfill$\square$\par\rm}

\def\boxing#1{\ \lower 3.5pt\vbox{\vskip 3.5pt\hrule \hbox{\strut\vrule
\ #1 \vrule} \hrule} }

\def\down#1{\ \lower 3.5pt\vbox{\vskip 3.5pt \hbox{\strut \ #1 \vrule} \hrule} }
\def\negdown#1{\ \lower 3.5pt\vbox{\vskip 3.5pt \hbox{\strut  \vrule \ #1 }\hrule} }

\def\items{\nl\leftskip = 25pt \parskip=\lineskip}
\def\enditems{\nl\leftskip = 0pt \parskip=\baselineskip}

\hsize=6.2 truein
\vsize=9 truein

\baselineskip=13 pt
\parskip=\baselineskip
 1

\parindent=0pt

\def\oover#1{\vbox{\ialign{##\crcr
{\npt o}\crcr\noalign{\kern 1pt\nointerlineskip}
$\hfil\displaystyle{#1}\hfil$\crcr}}}
\def\today{\number\day\space\ifcase\month\or January\or February\or
March\or April\or May\or June\or July\or August\or September\or
October\or November\or December\fi\space\number\year}

\newif \iftitlepage \titlepagetrue

\def\picture#1{\global\advance\diagramnumber by 1
$\epsfbox{doodlecatfig.\number\diagramnumber} \atop \hbox{#1} $}

\def\diagram#1{\global\advance\diagramnumber by 1
$$\epsfbox{doodlecatfig.\number\diagramnumber}$$\postdisplaypenalty=10000
\centerline{\npt\bf Figure\ \the\secnum.\the\diagramnumber \npt\sl\ #1}\par}

\newcount\diagramnumber
\diagramnumber=0

\newcount\chapternumber \newcount\diagramnumber 
\newcount\questionnumber \newcount\chapternumber
\chapternumber=1 \diagramnumber=0 \questionnumber=0
\newcount\secnum \secnum=0\newcount\subsecnum \subsecnum=0
\newcount\defnum \defnum=1\newcount\theonum \theonum=0
\newcount\lemnum \lemnum=0

\newcount\subsecnum
\newcount\defnum
\def\section#1{
                \vskip 10 pt
                \advance\secnum by 1  \subsecnum=0
                \leftline{\secfont \the\secnum \rm\quad\bf #1}
                }

\def\subsection#1{
                \vskip 10 pt
                \advance\subsecnum by 1 
                 \leftline{\secfont \the\secnum.\the\subsecnum\ \rm\quad\bf \ #1}
                }

\def\definition#1{
                \advance\defnum by 1 
               \par \bf Definition
\the\secnum .\the\defnum  \it \ #1\rm \par
                }
   
    \def\theorem#1{
                \advance\theonum by 1 
                \par\bf Theorem  \the\secnum
.\the\theonum \sl\ #1\rm\ \par
               }

\def\lemma#1{
                \advance\lemnum by 1 
                \par\bf Lemma  \the\secnum
.\the\lemnum \sl\ #1\rm \ \par
                }

\def\corollary#1{
                \advance\defnum by 1 
                \par\bf Corollary  \the\secnum
.\the\defnum \sl\ #1\rm\ \par
               }

\def\exercise#1{
                \advance\defnum by 1 
                \par\bf Exercise  \the\secnum
.\the\defnum \sl\ #1  \par\rm
               }

\def\cite#1{{\secfont [#1]}\nobreak}

\def\pagebreak{\vfill\eject}

\centerline{\bf Version \today}
\vglue 20 pt
\centerline{\secfont Planar Doodles:}
\centerline{\secfont Their Properties, Codes and Classification}
\medskip

\centerline{\nam Andrew Bartholomew}
\centerline{\nam Roger Fenn}
\centerline{\aff School of Mathematical Sciences, University of Sussex}
\centerline{\aff Falmer, Brighton, BN1 9RH, England}
\centerline{\aff e-mail: rogerf@sussex.ac.uk}
\centerline{\aff e-mail: andrewb@layer8.co.uk}
\vglue 20 pt


\centerline{\nam ABSTRACT}
\leftskip=0.5 in
\rightskip=0.5in

We present those properties of planar doodles, especially when regarded as 4-valent graphs, that enable us to classify them into {\it prime} and {\it super prime} doodles by 
analogy to a knot sum.  We describe a method for partialy characterising a doodle diagram by a {\it doodle code} that describes the complementary regions of the diagram and use 
that code to enumerate all possible prime and super prime doodle diagrams via their dual graph.  In addition we explore the relationship between planar dodles and twin groups, 
and note that a theorem of Tutte means that super prime doodles have a Hamiltonian circuit.  We hope to expand upon this last point in a follow-up paper.

\leftskip=0 in
\rightskip=0 in
\baselineskip=13 pt
\parskip=\baselineskip

\parskip=\baselineskip
\smallskip

\medskip
{\bf Keywords:} {\it Planar doodles, generalised knots, connectivity conditions, dual graphs, hamiltonian circuits, planarity conditions, doodle tables.}

{\bf MSC:} 05C10, 05A05, 05B30, 57M15, 57M25.

\section{Introduction}

Doodles are represented by  a collection of tame, immersed, oriented circles on $\Sigma$, an oriented surface. The curves are in general position with no triple or higher multiple points. Alternatively they are represented by 4-valent graphs in $\Sigma$, provided we take note of the cyclic ordering at vertices and the orientation of the edges. 

If we dictate that all the complementary regions of the graph are polygons and have at least 3 boundary edges then this picture completely and uniquely defines a doodle. This is called the {\it minimal representation}. A minimal representation is called {\it reduced} if we then throw away certain unnecessary components.

This paper arose out of a programme to try and classify planar doodles with a small number of crossings using the fact that a dual graph exists and is a plane graph. The authors with N. Kamada and S. Kamada have shown how it is possible  to produce tables of doodles on arbitrary surfaces, \cite{BFKK2}. However these methods are profligate with the genus of the ambient surface. 

It turns out to be necessary that we concentrate the classification according to certain connection properties of the doodle. For example, doodles which are disconnected by the removal of one or two crossings are complicated by the {\it Kishino effect}, see \cite{FTur}. That is, two trivial doodles may be summed, by analogy with classical knots, to produce a non trivial doodle.

We therefore restrict our attention to {\it prime doodles} which are not disconnected by the removal of 2 crossings and {\it super prime doodles} which are not disconnected by the removal of 3 crossings. These doodles have interesting properties which are described in what follows. For example, prime and super prime doodles have regions which  meet any other region, if at all, in either one crossing or one edge. It also means that the dual graph of such a doodle may be obtained from the 1-skeleton of a cell decomposition of a disc whose 2-cells are four-sided.

Another interesting fact which follows from a theorem of Tutte, \cite{T}, is that super prime doodles have a hamiltonian circuit. This means they have a representation by a particularly simple code. Examples are given in the text.

A computer search to enumerate graphical decompositions of a disc and identify planar graphs with four-sided regions is used to draw the corresponding doodle.  Details of how this is done will be described in section 8. This leads to a table of all suitable doodles up to 14 crossings. 

The contents of this paper are
\items

\item{1.} Introduction
\item{2.} Concepts and definitions
\item{3.} Connections
\item{4.} Infinite regions
\item{5.} Twins and doodles
\item{6.} Doodle codes
\item{7.} Dual graphs and their properties
\item{8.} Computer search strategy
\item{9.} Table: prime and super prime doodles up to 14 crossings

\enditems

\section{Concepts and definitions}

We can regard a doodle diagram as a connected 4-valent plane graph in some surface $\Sigma$ or as a collection of tame immersed oriented circles on $\Sigma$ in general position with no triple or higher multiple points. It is not 
difficult to see that the two definitions are equivalent provided we take note of the cyclic information at each vertex of the 4-valent graph and if we ignore the orientation of the circles. In this paper we will restrict to {\it 
planar doodle diagrams} where $\Sigma$ is the sphere which we think of as the plane with a point at infinity.  We regard a simple closed curve as the trivial doodle; in what follows all doodles are non-trivial unless explicitly stated 
otherwise.

Each immersed circle is called a {\it component} of the diagram.  We will use the term {\it vertex} and {\it crossing} interchangeably in what follows.
We call the components of the complement of a doodle diagram the {\it regions} of the diagram.   If the doodle diagram is connected then it defines a cell subdivision of the sphere in which the 2-cells are the regions. If a region 
contains  contains $i$ edges in its boundary curve then we refer it as an {\it $i$-gon} region.  In particular we speak of monogons, bigons and trigons when referring to $1$-gon, $2$-gon and $3$-gon regions respectively. 
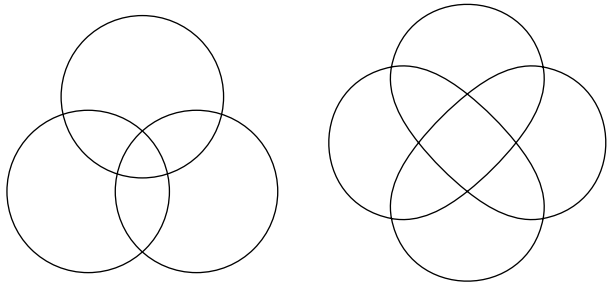
\diagram{Two doodle diagrams: the Borromean rings and the poppy}
Two example planar doodle diagrams are shown in Figure 2.1. The left hand one is called the {\it Borromean rings} which has three components and the right hand one with one component is called the {\it poppy}.
The regions of the Borromean rings are all trigons and there are eight of them. The poppy regions consist of eight trigons and two 4-gons.

\subsection{EIP, a simple but useful property}

Suppose we have a doodle diagram $D$ and a circle $C$ which is the boundary of a disk, $\Delta$ and $C$ meets $D$ is a finite number of points. Let $f:S^1\to S^2$ define a component of $D$ and pick $t\in S^1$ so that $f(t)\in C$. Let $N$ be a small neighbourhood of $t$ and suppose that $f(N)$ meets $S$ transversely. So $f(N)$ as an arrow either points into $\Delta$ or points out. 
We assume this for all such $t$ and every component. Then, an equal number of arrows point into and out of $\Delta$.

We call this simple but effective tool, the {\it Equal Intersection Property} or EIP for short. It will be used in the coming proofs.

\subsection{The doodle moves}

Two doodle diagrams are equivalent if they are related by a finite sequence of the moves $R_1, R_2$ shown in Figure 2.2 The equivalence class is called a doodle; however, following the usual convention, we will often not distinguish between a doodle and its diagram. 
\footnote{$ ^{(1)}$}{By theorem 2.1 a reduced doodle can be identified with its diagram}
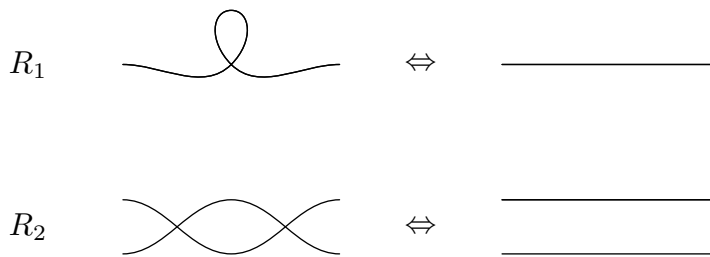
\diagram{The two allowable moves on a doodle diagram}
These two moves are analogous to the first two Reidemeister moves of knot theory but note that the third analogous move, illustrated below, is not allowed.
\diagram{Forbidden move}
The concept of a doodle is due to Fenn and Taylor \cite{FT}. Their definition differs from the one above, which was also used by Khovanov, \cite{Kh}, in that, in \cite{FT}, all component curves are simple, and so lack self crossings.  Doodles may also be defined on more general surfaces, whereupon they may be interpreted as virtual doodles, see \cite{BFKK1} for details but we shall not consider this more general doodle in this paper, however it is interesting to note that the complete 5-graph, $K_5$, exists as a doodle on a torus.

A circular component of a doodle which spans a disk whose interior contains no crossings is called a {\it floating circle}, which can be introduced or deleted at will. A doodle which has a disconnected diagram is called {\it reducible}. Otherwise the doodle is {\it irreducible}. We need only consider irreducible doodles. Clearly, irreducible doodles do not have floating circles.

The following theorem was proved in \cite{FT, Kh}.

\theorem{Any irreducible doodle has a unique diagram with the least number of crossings. Moreover this can be achieved by $R_1$ and $R_2$ moves which always decrease the number of crossings.\qed}
We call the unique diagram described by Theorem 2.1 {\it minimal}. So a { minimal doodle diagram} is defined to be one which is connected and with no monogon or  bigons.

A component of a minimal diagram is said to be {\it removable} if its absence leads to another minimal diagram, otherwise it is called {\it immovable}.
A component that is a small circle surrounding a crossing is called a {\it vertex circle}.  It is possible to add a vertex circle to any crossing of a minimal diagram and obtain another minimal diagram; a minimal 
diagram containing no removable vertex circles is called {\it reduced}. The smallest doodle containing a removable component that is not a vertex circle has twelve crossings.  Both the doodles $S12^4_1$ and $S12^4_2$ shown in section 9
are examples; these doodles may be constructed from the Borromean rings by adding another component that encircles more than one crossings.

\section{Connections}

Recall that a graph is called {\it $k$-connected} if it has more than $k$ vertices and is not  disconnected whenever fewer than $k$ vertices are removed. 
\footnote{$ ^1$}{Remember that connection in a graph is between two {\it vertices}. So the removal of the two end points of an edge need not disconnect a graph if the remaining vertices can be joined by a path.}
For example 1-connected is the same as connected and if a graph is 1-connected but not 2-connected then it has a {\it cut vertex} whose removal disconnects the connected graph. We extend these connected definitions to a general cell complex by saying that it is $k$-connected if its 1-skeleton is $k$-connected.

A doodle with a cut vertex can be thought of as the {\it one point union} of two doodles: the joining points of the doodles lying in the interior of an edge from each diagram. 

A doodle with no cut vertices is {\it 2-connected}. It is easy to see that for 2-connected diagrams all regions are embedded $i$-gons for some $i$, meaning that there are $i$ distinct edges and $i$ distinct vertices.
We will always assume unless expressly stated that for the rest of this paper all doodles are 2-connected. 

A 3-connected doodle is called {\it prime}. The examples in Figure 2.1 are 3-connected and hence prime. We consider prime doodles because our methods do not provide a way of enumerating doodles which are not prime. These are doodles which are disconnected by the removal of two vertices. By analogy with classical knot theory these appear to be the {\it sum} of two doodles.
Un-prime doodles are not uninteresting and an un-prime doodle need not determined by its summands.  The example in Figure 3.4 is 2-connected but not 3-connected,  and its summands are two trivial doodles. This phenomenon also occurs for virtual knots with the Kishino knot, see \cite{FTur}.
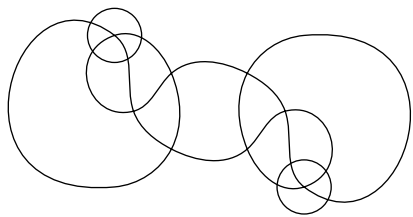
 \diagram{Sum of trivial doodles}
 A 4-connected doodle is called {\it super prime}. The doodles in Figure 2.1 are super prime but the doodles in Figure 4.6 are prime but not super prime. We cannot extend to any higher connectivity because each crossing has valency four.

\section{Infinite regions}

We may consider a spherical doodle to be planar by choosing an {\it infinite region} to contain the point at infinity and the other regions are then referred to as {\it finite} regions.  If the doodle is 2-connected then the complementary closure of the infinite region, i.e. the union of the finite regions, is a disk, $\Delta$. The boundary of this disk, $\partial \Delta$, will be the {\it boundary} of the doodle for this choice of infinite region.  

Call a finite region whose closure has a non-empty intersection with the boundary a {\it boundary region} of the doodle.  Consider the crossings on the boundary, they are incident with two edges forming part of the boundary and two edges which point into the interior of the boundary. There are locally four regions which meet a crossing $v$ on the boundary. In cyclic order there is the infinite region, a finite region which meets the boundary in an edge $e_1$ say, a finite region which meets the boundary in a vertex $v$ and a finite region which meets the boundary in an edge $e_2$. The edges $e_1$ and $e_2$ are boundary edges meeting at $v$.

\lemma{Let $D$ be an irreducible minimal  2-connected doodle with a chosen infinite region. Let $R$ be a boundary region of $D$. Then $R$ meets the boundary either in isolated boundary vertices or a number of disjoint edges but not both.}
{\bf  Proof:}  Suppose the conclusion of the lemma is false and $R$ meets the boundary in the isolated vertex $v$ and a disjoint edge $e$. Then there is a circle $C$ which passes through $v$ and transversely through an interior point of $e$ but is otherwise disjoint from $D$. But this means that there are three edges pointing into the interior of $C$ which is an impossibility by the EIP. \qed

The finite regions can be subdivided into three kinds by the above lemma.
We say that a finite region, $R$, of a minimal prime doodle diagram $D$ is a {\it vertex boundary region} or {\it v-region} if it meets the boundary of $\Delta$ in vertices, an {\it edge boundary region} or {\it e-region} if it meets the boundary of $\Delta$ in edges and an {\it interior region} otherwise. The v-regions and e-regions alternate around the boundary of $\Delta$. If the infinite region is a $p$-gon then there are $p$ edges and $p$ vertices on the boundary of $\Delta$. Every boundary vertex is connected to two boundary edges and two edges which point into the interior of $\Delta$. These last two edges are called {\it semi boundary edges}. Each semi boundary edge is in the boundary of a $v$-region and an $e$-region. 

\lemma{Let $D$ be an irreducible minimal  2-connected doodle with a chosen infinite region. If two boundary regions meet in an edge, then one region is a v-region and the other is an e-region}
{\bf  Proof:} If the regions meet in a semi boundary edge then they are alternate boundary regions. If they meet in an interior edge and they are both e-regions then there is a simple closed curve meeting the doodle in three points contradicting the EIP. A similar argument holds if both are v-regions.
\qed

\lemma{The boundary regions of a prime doodle either meet the boundary in one edge or one vertex.}
{\bf Proof:} We have already seen from lemma 4.1 that a boundary region meets the boundary either in edges or vertices. If there are two edges in the region meeting the boundary then there is clearly a pair of vertices which disconnect the doodle contrary to the fact that the doodle is prime and hence 3-connected. A similar argument works for the case where the boundary regions meets the boundary in crossing points. \qed

\theorem{A doodle is prime if and only if every pair of  regions  are either disjoint, meet in a vertex or meet in an edge.}
{\bf Proof:} Since any region can be chosen as the infinite region we have the implication one way from lemma 4.3. 

Suppose, conversely, that the doodle is not prime and is 2-connected but not 3-connected. Then there are two crossings, $A, B$, whose removal disconnects the doodle into two parts $D_1$ and $D_2$. 
Assume we have chosen an infinite region and $A$ lies in the boundary. Because $A, B$ disconnects the doodle, $B$ must also lie in the boundary.
There is a simple closed curve $C$ which only meets the doodle in $A, B$ and such that $D_1$ say, lies inside $C$ and $D_2$ lies outside. There are two cases, illustrated in Figure 4.5, depending on whether,
at each of $A$ and $B$, an odd number of edges or an even number of edges point inside $C$ which is drawn in red.

\diagram{2-connected but not 3-connected}

The sub-arc of $C$ joining $A$ to $B$ which lies inside the doodle lies in a region which meets the infinite region on either two edges or two crossings. 
\qed
\corollary{Every semi boundary edge of a prime doodle connects a boundary vertex to an interior vertex.}
{\bf Proof:} If  both vertices were boundary vertices then their removal would disconnect the doodle contrary to primeness. \qed

\subsection{Inner complements}

If $R$ is a region of a prime doodle chosen to be infinite then its {\it inner complement}, $R^c$, is defined to be all the interior vertices, edges and regions which have no contact with the boundary of $R$. Examples of inner complements in red and boundary edges in blue, are given in the Figure 4.6. The semi-boundary arcs are black and join the red crossings to the blue crossings.
\diagram{Inner complements} 

\lemma{Let $D$ be a prime doodle and let $e_1, e_2$ be semi-boundary edges which meet in an interior vertex $v$ and are edges of a region $S$. Then $S$ is an $e$-trigon.}
{\bf Proof:}  The region $S$ is a boundary region. If $S$ is a $v$-region and meets the boundary in $v'$ then $v, v'$ are the vertices of a bigon contrary to minimality. So $S$ is an $e$-region and meets the boundary in an edge which forms a trigon with $e_1, e_2$. \qed

\lemma{Let $D$ be a prime doodle. With the notation above $R^c$ cannot have a vertex of valency 0 or 1.}
{\bf Proof:} The vertex if it exists will be the common vertex of at least three semi-boundary edges. By lemma 4.4 one of the semi-boundary edges must be in the boundary of two adjacent $e$-trigons. But this is impossible as $e$-regions and $v$-regions alternate. \qed

\theorem{The inner complement of any region of a prime doodle is acyclic ($H_1=0$) and contains at least one region.}
{\bf Proof:} Let $R$ be the infinite region. Then $R^c$ must be non-empty because the doodle is prime. By lemma 4.5 any interior vertex must meet at least two edges of $R^c$.
This means that any path component of the  1-skeleton of $R^c$ cannot be a tree as this has a vertex of valency one. So the 1-skeleton contains a cycle whose interior contains at least one interior region. Similarly there can be no holes in $R^c$ because they are filled with interior regions. \qed

The existence of one interior region is best possible. There are infinite sequences of doodles with one interior region.

\theorem{The inner complement of any region of a super prime doodle is a disk}
{\bf Proof:} If the complement has a cut vertex then the doodle can be disconnected with three vertices. So the complement is the disjoint union of disks. But the complement must be connected by $2$-connectivity. \qed
We now have a picture of the body of a super prime doodle as a disk with an inner complement a disk separated by an annular region.

\subsection{Hamiltonian circuits}

Recall that a topologically embedded circle in a graph which contains all the vertices of the graph is called a {\it hamiltonian circuit}. A theorem of Tutte, \cite{T}, states that all 4-connected planar graphs contain a hamiltonian circuit. This means that all super prime doodles have a hamiltonian circuit.  Whilst many prime doodles have a Hamiltonian circuit, it is an open question whether this is true of all prime doodles.

If a doodle has a hamiltonian circuit then we can give the following description and associated code in the form of cycles. Assume that the doodle has $n$ crossings and $2n$ edges. Imagine that the hamiltonian circuit is the equator of a sphere. Then the $n$ edges which are not in the hamiltonian circuit either cross the northern hemisphere or the southern hemisphere. 

In figure 4.7 we look at the Borromean rings, the poppy, the nine crossing doodle $P9^1_1$ and the ten crossing doodle $S10^1_1$ We view them as if we are looking down from the North pole. The edges in the northern hemisphere are coloured blue and the edges in the southern hemisphere are coloured red.

The associated cycles are defined as follows. Label the equatorial crossings $0, 1, \ldots, n-1$ and write $i+j$ if $i$ is joined to $j$ by a red arc and $i-j$  if $i$ is joined to $j$ by a blue arc, (all mod $n$). 
Starting from 1 we follow the arcs until we return to 1. If any crossings are unvisited we choose the smallest and continue again as if we were writing a permutation as a product of disjoint cycles.
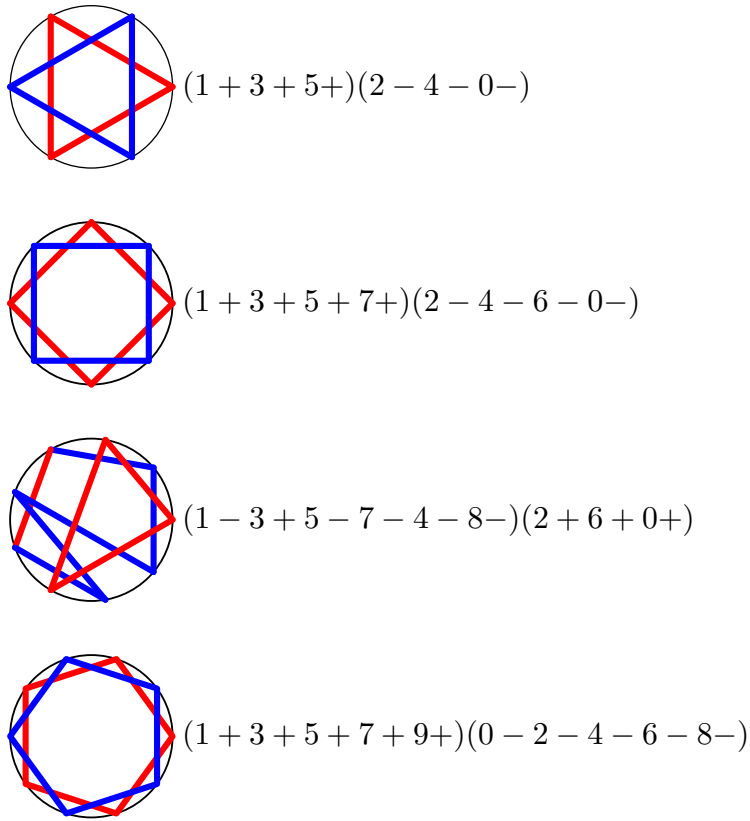
\diagram{4 doodles with hamiltonian circuits and associated cycles}

We hope to return to the subject of Hamiltonian circuits in a future paper, as outlined in section 10.

\section{Twins and doodles}

As the title suggests, in this section we will define the twin groups and show how their closures can make a diagram which represents any doodle. We will sketch the main points which are easily understood but full details can be found in the following references, \cite{BF} 

\subsection{The twin group $TW_n$}

We can define the twin group as a presentation thus
$$TW_n=\langle t_1,\ldots, t_{n-1}\mid t_i^2=1, \hbox{ for all }i,\ t_it_j=t_jt_i, \hbox{ if } |i-j|>1\rangle
$$
This is a right angled Coxeter group.

Any particular word in the generators can be represented, braid style, as in Figure 5.7. The figure offers two solutions, the left hand side fits in with the decomposition of a doodle into a Seifert graph with horizontal lines or {\it threads} and the right hand side mirrors a more conventional braid style with the {\it strings} picked out in various colours.
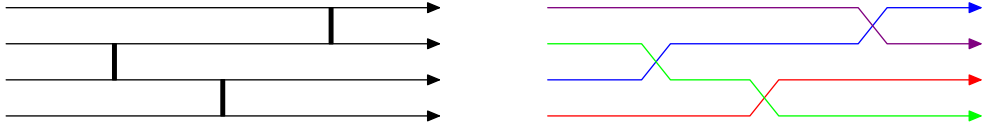
\diagram{Representing the word $t_2t_1t_3$}
These diagrams can be closed as in the usual way for classical braids. The resulting doodle diagram is called {\it braided}. Any doodle can be represented by a braided diagram, see \cite{F} For example the Borromean rings are the closure of $(t_1t_2)^3$ and the poppy is the closure of $(t_1t_2)^4$

However a minimal diagram may not be braided and the following figures will illustrate this phenomenon and how it is resolved.

The doodle diagram illustrated  in Figure 5.8 is minimal but is not braided.
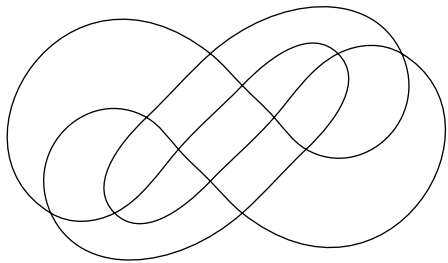
\diagram{Minimal diagram but not braided}
Although it is not braided it nevertheless has an interior which is a disk.

We now orient the doodle and convert it into a {\it Seifert graph} by the the following move which converts crossings into bridges, (red in the figure).
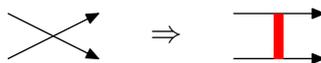
\diagram{The Seifert move}
If we now apply this move to the crossings of a doodle diagram the diagram becomes a union of disjoint simple cycles joined by bridges with the orientation respected as in Figure 5.9.
The result for the above doodle looks like Figure 5.10.
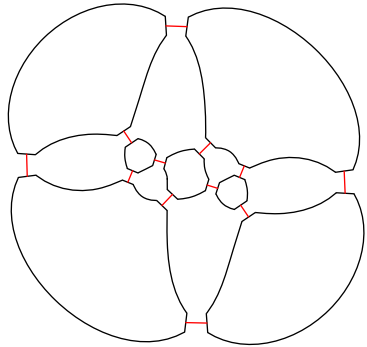
\diagram{Seifert graph with bridges in red}
In order to convert this picture to the closure of a twin it is necessary to apply yet another move, the {\it V-move} illustrated in Figure 5.11
\diagram{The V-move}
The two arcs to the left of Figure 5.11 are from {\it different} cycles and are oriented in different directions. After the move the two cycles become one and a new cycle is born connected to the amalgamated cycles by two bridges. When the V-move is applied in all situations the cycles become nested, (remember we are on a sphere) and the twin word, defined up to a cycle, can be read off.

Applying this to our example we have the word
$w=t_1t_2t_3t_4t_3t_2t_3t_2t_1t_2t_3t_4t_3t_2t_3t_2$
 and the corresponding diagram in Figure 5.12 representing its closure.
\diagram{The closure of $w$}

\section{Doodle codes}

Suppose a minimal doodle diagram has $n$ crossings and $f_i$ $i$-gon regions, for $3\le i\le p$.  Then, the diagram has $2n$ edges, and by Euler's formula the number of regions is $f=n+2$, one of which is infinite. So\nl
{a)} $f_3+f_4+\cdots +f_p= n+2$ \nl
{b)} $ 3f_3+4f_4+\cdots+pf_p=4n$ \nl
Conversely, a doodle diagram is partly characterised by a solution to a) and b) for a given number of crossings.  We call the combination of the number of crossings
and a solution to a) and b) a {\it doodle code}.

For example, the codes of the Borromean rings and the poppy are $6,8$, six crossings with eight trigons, and $8,8,2$, eight crossings, eight trigons and two $4$-gons, respectively.

The following two lemmas place constraints on the number 

\lemma{The number of trigons of a minimal doodle with $n$ crossings satisfies
$$8\le f_3\le{4n\over 3}$$
and the number of $i$-gons where $i>3$ satisfy
$$0\le f_i\le {{4n-24}\over i},\ i>3$$}
{\bf Proof}  If we multiply equation a) above by 4 and subtract b) we see that 
$$f_3=8+f_5+2f_6+\cdots+(p-4)f_p$$
where $p$ is the size of the largest polygon region and the left inequality of the first condition follows. The right inequality of the first condition follows from b).   The proof of the second inequalities follows similarly. \qed

\lemma{If a doodle is prime with $n$ crossings then every region has less than $(n+1)/2$ edges.}
{\bf Proof}  Using the above notation, suppose that the infinite region has $p$ edges and $p$ vertices. Then this defines $2p$ boundary regions which are all distinct because the doodle is prime and so no boundary region can have 2 edges (vertices) in common with the infinite region.
So there are $n+1-2p$ interior regions. This is positive by theorem 4.3. \qed

\section{Dual graphs and their properties}

The dual graph, $\Gamma=\Gamma(D)$, of a doodle diagram $D$ is defined by associating a vertex of $\Gamma$ with each region of the $D$ and an edge of $\Gamma$ with each edge of $D$.  If an edge $e$ of $D$
lies between regions $R_1$ and $R_2$ then the corresponding edge $e'$ of $\Gamma$ joins the vertices corresponding to $R_1$ and $R_2$.
Since the doodle is planar, so is its dual graph; moreover, since $D$ is a four-valent graph, $\Gamma$ divides the sphere into 4 sided cells.

Assume that an infinite region of an irreducible prime doodle $D$ is chosen. The vertex of $\Gamma$ corresponding to the infinite region is called the {\it infinite vertex} and labelled $v_{\infty}$.  A vertex of $\Gamma$ corresponding to an e-boundary region of $D$ is called an {\it e-boundary vertex}. These vertices are joined to $v_{\infty}$ by edges of $\Gamma$, called {\it outside} edges,  forming the star sub-graph of $v_{\infty}$. 

A vertex of $\Gamma$ that corresponds to a v-boundary region of $D$ is called a {\it v-boundary vertex}. A cell with $v_{\infty}$ in its boundary is called an {\it infinite cell}. Besides $v_{\infty}$ it also has a v-boundary vertex separating two e-boundary vertices in its boundary. If the infinite region of $D$ is a $p$-gon then the link sub-graph 
of $v_{\infty}$, called the {\it boundary} of $\Gamma$ relative to $v_{\infty}$, consists of $2p$ edges alternatively joining v-boundary and e-boundary vertices.

The graph $\Gamma^* = \Gamma - v_{\infty}$ is called the {\it complementary dual graph} relative to $v_{\infty}$.
The outside edges point into the infinite region of the plane and are incident with only one vertex other than $v_{\infty}$, which is an e-boundary vertex. A vertex of $\Gamma$ which is not $v_{\infty}$, an e-boundary vertex or a v-boundary vertex is called an {\it interior vertex},

Figure 7.14 shows the complementary dual graph of the poppy; it has four infinite cells, four outside edges, four e-boundary vertices, four v-boundary vertices and one interior vertex.

\diagram{complementary dual graph $\Gamma^*$ of the poppy}

\lemma{A doodle is prime if, and only if, for every choice of infinite region, the boundary of the dual graph is an embedded circle.}
{\bf Proof}  If a doodle $D$ is not prime then, as in figure 4.5, the doodle has two disconnecting vertices.
Consequently the complementary dual graph $\Gamma^*$ of $D$ has a cut point as shown in red in Figure 7.15. So the boundary of the dual graph cannot be an embedded circle.
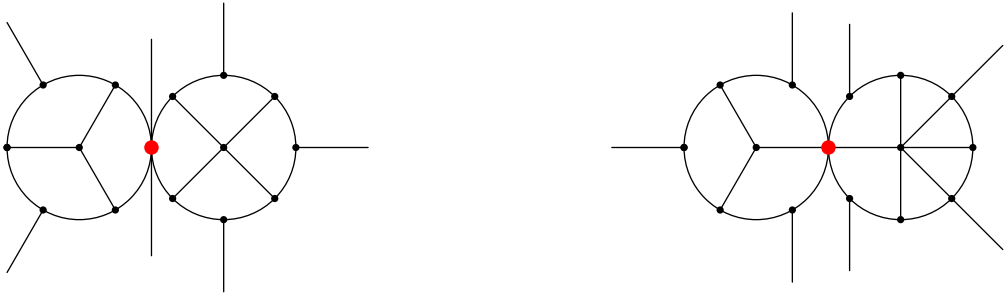
 \diagram{Dual graphs with cut points}

For the converse suppose there is a choice of infinite region for which the boundary of the dual graph, $\partial \Gamma$, is not an embedded circle.  Since $D$ is 2-connected, $\partial \Gamma$ must be connected, since 
considering the infinite cells in sequence around $v_\infty$ we may concatenate the edges in cell boundary that are not outside edges to form a boundary path that visits every vertex in $\partial \Gamma$. If $\partial \Gamma$ is not an 
embedded circle then this boundary path is not a simple closed curve, each self intersection occurring at either an e-boundary or a v-boundary vertex.  If an intersection point, $v$, is a e-boundary vertex then 
there is a region of $D$ that meets the infinite region in two distinct edges, contradicting Theorem 4.2.  If $v$ is a v-boundary 
vertex, there are two distinct infinite cells that both contain $v$ in their boundary and $D$ has a region that meets the infinite region in two distinct vertices, again contradicting Theorem 4.2.\qed

Let $D$ be a prime minimal planar doodle with $n$ crossings. Chose an infinite region with $p$ edges.  If $\Gamma$ is the dual graph of $D$ let $v_{\infty}$ be the vertex of  $\Gamma$ corresponding to the infinite region. Put
$\Gamma^* = \Gamma - v_{\infty}$,  the complementary dual graph, and $\Gamma' = \Gamma^* - \hbox{star}(v_{\infty},\Gamma)$.  By Lemma 7.8, $\Gamma'$ is the $1$-skeleton of a 
cellular decomposition of disc $\Delta$.
We may summarise the properties of the dual graph of a prime minimal doodle in the following lemma. 

\lemma{ Using the above notation:\nl
{1.} $\Gamma$ has $n+2$ vertices, $2n$ edges and $n$ regions\nl
{2.} There are $n+1$ vertices in $\Gamma'$\nl
{3.} There are $2p$ edges and $2p$ vertices in the boundary of $\Delta$\nl
{4.} There are $2n-3p$ interior edges and $n+1-2p$ interior vertices in $\Delta$\nl
{5.} There are $n-p$ regions in $\Delta$, all four-sided\nl
{6.} There are $f_p-1$ vertices of valency $p$ in $\Gamma'$ and for $i\ge3,\ i\ne p$, there are $f_i$ vertices of valency $i$ in $\Gamma'$\nl
{7.} Every region in $\Delta$ is adjacent to at least two other regions.\nl
{8.} There is at most one edge joining any pair of vertices of $\Gamma$.\nl
}
{\bf Proof:}  
Points 1. to 6. are easy consequences of previous discussions.\nl
For point 7., if there were a region, $\rho$, of $\Delta$ incident with just one other region then $\rho$ would have three 
edges in the boundary of $ \Delta$.  Since in $\Gamma$ alternate vertices are joined to the infinite vertex, there would then be a vertex of valency two in $\Gamma$.  Such a 
vertex would correspond to a bigon in the doodle, contradicting the fact that $D$ is minimal.  Point 8 follows immediately from theorem 4.2.\qed

\subsection{Constructing doodles from dual graphs}

Let $\Delta$ be a disc with $2p$ vertices in its boundary and $k$ vertices in its interior and construct a graph $\Gamma$ from $\Delta$ as follows.  On each alternate boundary vertex, construct an outward edge to a vertex, $v_\infty$, considered to be at the point of infinity.

Let $n=2p+k-1$ and suppose it is possible to add $2n-3p$ additional edges to the interior of $\Delta$ so that the resultant graph, $\Gamma$, is the $1$-skeleton of a cellular decomposition of the sphere whose regions 
all have four edges.  Then, we may recover a doodle $D$ with $n$ crossings whose dual graph is $\Gamma$ by constructing a crossing in each region of $\Gamma$, joining the midpoint of each edge in the boundary of a 
region, $\rho$, to the midpoint of the opposite edge; that is, of the unique non-adjacent edge in the boundary of $\rho$.

In Figure 7.16 are two examples of dual graphs and in Figure 7.17 a doodle is constructed from the left hand example.
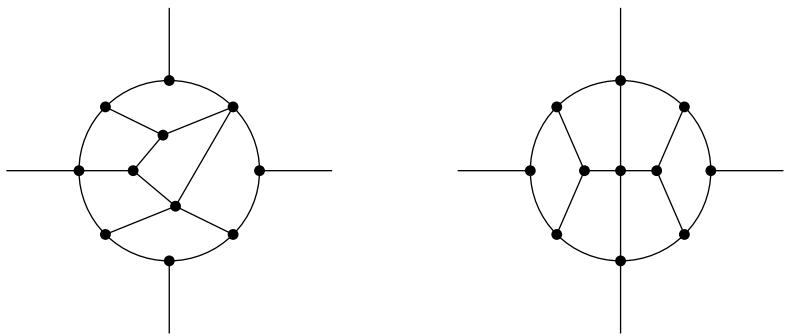
\diagram{Dual graph examples}
\diagram{Constructed doodle}

If the graph $\Gamma$ contains no vertex of valency two then $D$ is minimal.  Note that, although the boundary of $\Gamma$ is an embedded circle, without placing constraints 
on the placement of the $2n-4p$ additional interior edges, $D$ may not be prime since a different choice of infinite region may result in a boundary that is not 
an embedded circle.  

\section{Computer search strategy}

A computer search has been implemented to enumerate prime minimal doodles with up to fourteen crossings.  The methods used are general, 
although the processing time required above fourteen crossings becomes prohibitive.

The search begins by evaluating doodle codes.  By lemma 6.7, it is an easy matter to determine all possible values of $f_i$, the number of $i$-gons of prime doodles with a given number of crossings.  The following  table shows codes up to $n=11$.

$$\vbox{ \offinterlineskip \halign{\strut \vrule \hfil \quad #\quad \hfil \vrule
                                        & \hfil \quad #\quad \hfil \vrule 
                                        & \hfil \quad #\quad \hfil \vrule
                                        & \hfil \quad #\quad \hfil \vrule 
                                        \cr
\noalign{\hrule}
number of & \multispan3 number of i-gons \vrule \cr
\omit & \multispan3\hrulefill \cr
crossings & $f_3$ & $f_4$ & $f_5$ \cr
\noalign {\hrule}
$6$ & $8$ & $ $   & $ $ \cr
\noalign {\hrule}
$8$ & $8$ & $2 $  & $ $ \cr
\noalign {\hrule}
$9$ & $8$ & $ 3$  & $ $ \cr
\noalign {\hrule}
$10$ & $8$ & $4$  & $ $ \cr
\noalign {\hrule}
$10$ & $9$ & $2$  & $1$ \cr
\noalign {\hrule}
$10$ & $10$ & $0$ & $2$ \cr
\noalign {\hrule}
$11$ & $8$ & $5$  & $ $ \cr
\noalign {\hrule}
$11$ & $9$ & $3$  & $1$ \cr
\noalign {\hrule}
$11$ & $10$ & $1$ & $2$ \cr
\noalign {\hrule}}}$$

The search algorithm uses the construction in section 7.1 to create planar dual graphs with four-sided regions from each doodle code.  If the largest index $i$  for which $f_i \ne 0$ 
is $p$, the algorithm considers a disc $\Delta$ with $2p$ vertices in its boundary and $n+1-2p$ interior vertices. 
Since the algorithm searches for minimal doodles, each vertex in the dual graph must have valency at least three.  Every alternate vertex in the boundary of $\Delta$ is an e-boundary vertex 
and is joined by an edge to the infinite vertex, so the search considers all possible decompositions of the disc into four-sided regions such that every interior vertex
and every alternate boundary vertex has valency at least three.

Doodles are identified up to reflection in a line in the plane so the dual graph need only be constructed up to symmetry described by the dihedral group on $2p$ elements.

A dual graph $\Gamma$ is represented as a sparse matrix, $I$, that records the edges of $\Gamma \cap \Delta$ but not the edges to $v_{\infty}$.  Thus $I$ is an 
$n+1 \times n+1$ matrix which, by lemma 7.9 8.\ has entries equal to $1$ or $0$.  A $1$ in entry $(r,c)$ indicates that vertices $r$ and $c$ are joined by an edge, so $I$ is a self transpose matrix 
with zeros along the leading diagonal.  The vertices in the boundary of $\Delta$ are numbered $0,\ldots,2p-1$ clockwise around $\Delta$ and we assume, without loss in generality, that vertex $0$
is a v-boundary vertex.  Interior vertices are numbered as described below.  Clearly, a permutation of the vertex numbering does not change the doodle determined by an incidence matrix.

The matrix $I$ is initialised with the edges in the boundary of $\Delta$.  In order to reduce the number of matrices to be considered to a managable level, $I$ is constructed in stages as the 
remaining $2n-3p$ edges are added.  The stages are determined by the combinatorial characteristices of the dual graph, as described in the following sections.

\pagebreak
\subsection{Stage 1}

As noted above, an e-boundary vertex and is joined by an edge to the infinite vertex so will therefore have valency at least three.  
The v-boundary boundary vertices must be joined either to an interior vertex or to another boundary vertex.  An edge joining a v-boundary 
vertex to another boundary vertex is called a {\it boundary chord}.  A boundary chord divides $\Delta$ into two regions whose closure must each contain both endpoints of any edge 
of $\Gamma$ they intersect, since $\Gamma$ is required to be planar.

The first stage of the enumeration is to identify all possible configurations of boundary chords, including the case where there are none, and to divide $\Delta$ into {\it chord regions} that 
have the property that a vertex in the chord region's boundary only connects to an interior vertex.  Note that, by lemma 7.9 8, we do not need to consider boundary chords between adjacent vertices.

For the case where there are no boundary chords we pass a copy of $I$ containing entries corresponding to the boundary edges only to the next stage of the enumeration.
The other cases require a systematic enumeration of boundary chord patterns.

We represent a boundary chord by a pair of integers $(a,b)$ with the property that $a<b$, indicating that the corresponding boundary vertices are joined by the chord.  Boundary chords are ordered by definining 
$(a,b) < (c,d)$ if, and only if, $a<c$, or $a = c$ and $b < d$.  A pattern, $P_k$, of $k$ boundary chords is an ordered sequence of chords 
$\{(a_i,b_i): (a_i,b_i) < (a_j,b_j), i,j=1,\ldots,k, i<j\}$. 
Patterns of boundary chords are ordered by defining $P_k = \{(a^k_i,b^k_i)\} < P_m = \{(a^m_i,b^m_i)\}$ if, and only if, $k<m$, or $k=m$ and $(a^k_i,b^k_i) < (a^m_i,b^m_i)$ for some $i$.  Clearly,
the boundary chords in a pattern may intersect only at their endpoints.

\lemma {A boundary chord in the dual graph of a planar doodle $D$ joins an e-boundary bertex to a v-boundary vertex.}
{\bf Proof} The presence of a boundary chord in the dual graph of $D$ means that there is a simple closed curve $C$ that meets $D$ in exactly three points: the edge corresponding to the chord and an 
edge or a crossing corresponding to the boundary vertices it joins, depending on whether it is an e-boundary vertex or a v-boundary vertex respectively.  Thus, if the chord joined two e-boundary vertices or two
v-boundary vertices it would violate the EIP. \qed

A chord region may be regarded as an i-gon whose vertices are all boundary vertices and whose edges are either boundary edges of $\Delta$ or boundary chords.  The following lemma restricts the boundary chord patterns
that need to be considered.

\lemma {If an i-gon is has a cell decomposition where all the 2-cells are four-sided regions then i must be even.}
{\bf Proof} If the i-gon is divided into $r$ four-sided regions and the 1-seleton of the cell decomposition has $e$ edges then $4r = 2e - i$, so $i$ must be even. \qed

By lemma 8.11 and the dihedral symmetry of dual graphs, we need only consider patterns comprised of a single boundary chord of the form $\{(0,b)\}$, where $b=2i+1, i=1,\ldots$ and $b<p$, since we may 
rotate and reflect any other single-chord pattern capable of being divided into four-sided regions into a pattern of this form. 

The enumeration of boundary chord patterns is constrained by the minimum number of interior vertices required to divide the chord regions into four-sided cells.  If a chord region has six or more edges in its 
boundary then it is possible to divide it into four-sided regions by adding a single interior vertex, since each alternate boundary vertex may be joined to the interior vertex, which will have valency at 
least three. If the chord region has four boundary edges we use the following lemma.

\lemma {If a 4-gon is to be divided into four-sided regions by adding $k>0$ interior vertices such that a) every interior vertex has valency at least three and 
b) no interior vertex is joined to another vertex by more than one edge, then $k \ge 4$.}
{\bf Proof} Suppose the vertices of the 4-gon are $a,b,c,d$ and that $a$ is joined to an interior vertex $e$.  By lemma 8.11, vertex $e$ cannot be joined to $b$ or $d$.  Since $e$ is required to have valency at least three, 
if $e$ is not joined to $c$ there must be two other interior vertices $f,g$ to which it is joined.  By lemma 8.11 $f,g$ cannot be joined to each other, nor to $a$ or $c$, and since they cannot both be joined to 
both $b$ and $d$, there must be a fourth interior vertex $g$.  If $e$ is joined to $c$ then it must be joined to another disjoint interior vertex $f$ and we may recursively repeat the above argument with 4-gon 
$a,e,c,b$ or $a,e,c,d$, whichever contains $f$ and the result follows. \qed

Adding a boundary chord to a chord pattern increases the number of interior vertices required by at least one.  Therefore, the enumeration begins with a list of single chord patterns of the above form that do not
require too many interior vertices.  Iteratively, the first pattern, $P_k$, on the list is extended by considering the ordered list of boundary chords larger than the last chord in $P_k$.  Sequentially, each such 
chord is added to $P_k$ to form a pattern $P_{k+1}$ and the minimum number of interior vertices required to divide the chord regions of $P_{k+1}$ into four-sided cells is evaluated. If $P_{k+1}$ does
not require too many interior vertices, it is added to the end of the list.  Finally, the edges corresponding to $P_k$ are recorded in the first $2p\times 2p$ submatrix of a copy of $I$, which is then passed to the next 
stage of the enumeration and the pattern $P_k$ removed from the list.

\subsection{Stage 2}

The boundary of a chord region may contain v-boundary vertices that are not the endpoint of a boundary chord.  Such v-boundary vertices are called {\it inside vertices} becasue they must 
be joined to an interior vertex to have a valencey of at least three.  
The second stage takes a combination of zero or more boundary chords and distributes the $n+1-2p$ interior vertices across the chord regions in all possible ways.  For each allocation of interior 
vertices the chord regions are considered, first in the order in which they are are incident with the edges in the boundary of $\Delta$ ordered clockwise from vertex $0$.  Then, if there are chord regions that 
only meet the boundary of $\Delta$ in vertices, in the order in which they are incident with the vertices in the boundary of $\Delta$, ordered clockwise from vertex $1$.  
For each region in turn, a choice of interior vertex made for connecting each inside vertex.  This choice is called the {\it interior choice} for the inside vertex.

The possible interior choices are enumerated, initialised by attaching as many inside vertices to the same interior vertex as possible.  The inside vertices of a chord region are  
boundary vertices, which are numbered as above, and the order in which the interior choices are made respects that order.  By lemma 7.9 6.\ the doodle code determines the maximum 
valency of a vertex in the dual graph, which determines how many inside vertices can have the same interior choice.  The process of making an 
interior choice for the inside vertices in the above order implicitly numbers the interior vertices that are required for the choice.  If, given a set of interior choices, there are two or more unattached interior 
vertices and we are going to advance the enumeration so that an unattached vertex becomes joined to the boundary, then we do not need to select each available unattached interior vertex 
in turn, since doing so effectively renumbers the unattached interior vertices.  At each step of the enumeration the interior choice edges are recorded in the right-most $n+1-2p$ columns of the first
$2p$ rows of a copy of $I$ which is then passed to the next stage of the enumeration.

\subsection{Stage 3}

The remaining edges in an incidence matrix produced by the first two stages either join a boundary vertex to an interior vertex or join two interior vertices.  Since the dual graph has at most one edge 
between a pair of vertices, it is straightforward to identify the zero locations of $I$ in the right-most $n+1-2p$ columns that lie to the right of the leading diagonal and then select the appropriate number 
of locations to set to 1 in all possible ways.

The third stage of the enumeration iterates through the possibilities, rejecting any that result in the creation of triangular regions: that is, a location in $I$ that corresponds to an edge between 
vertices $u$ and $v$ is not selected if $u$ and $v$ are already joined to a common vertex.  Similarly, combinations of locations are ignored if they result in a vertex having valency less than three.

Each completed matrix $I$ is checked to ensure the valency requirements of the doodle code are met and that the corresponding graph $\Gamma$ is planar. For a planar graph, we may trace the 
boundary of a region by moving along an edge towards a crossing and turning right, since there is a unique choice for the subsequent edge.  In this manner, the regions of planar graphs are checked to 
ensure they are all four-sided and, if so, we say that the incidence matrix is {\it admissible}.

\lemma{The dual graph represented by an admissible incidence matrix $I$ is connected.}

{\bf Proof} Each component of the dual graph represented by an admissible matrix $I$ determines a cellular decomposition of $S^2$ with four-sided regions.  Let the number of
vertices, edges and faces of the decomposition determined by the $k$-th component be $v_k$, $e_k$, and $f_k$ respectively.  Then, $4f_k = 2e_k$ and $v_k-e_k+f_k = 2$, so
$f_k = e_k/2$ and $v_k = 2 + e_k/2$.  In total the dual graph has $n+2$ vertices and $2n$ edges, so $n+2 = \Sigma_k v_k = 2k +1/2\Sigma_k e_k = 2k +n$ and therefore $k=1$. \qed

By lemma 8.13 the dual doodle of an admissible incidence matrix is irreducible.  The doodle is minimal since the matrix defines a graph whose vertices all have valency at least three.
Since some doodles have multiple regions with the most number of edges, the algorithm may produce different incidence matrices for the same doodle, corresponding to the different diagrams
obtained by a different choice of infinite region.  There is also the possiblility of duplication due to vertex renumbering as the remaining edges are added, if there are multiple unattached
interior vertices in a chord region.
The unoriented left preferred Gauss code for a doodle defined in \cite{BFKK2} provides a unique representation for irreducible minimal doodles that allows a list of distinct dual doodles to be determined, which may be classified as either prime or super-prime using theorem 4.4.

Some of the doodles produced by the search include removable vertex circles and have been discarded, so that the search yields only reduced minimal primes.
These doodles have been identified manually by inspecting diagrams obtained from their Gauss code, using an algorithm based on Thurston's circle packing algorithm, as implemented by Knotscape \cite{HT} and the first author's draw programme \cite{B}.

As noted in section 7.1, it is possible to construct a dual graph whose boundary is an embedded circle for an non-prime minimal doodle, so it is possible the
algorithm produces an non-prime doodle.  However, in light of lemma 7.8, the construction is guaranteed to find all those doodles that are prime.

All of the doodles produced up to fourteen crossings are prime, so it is an interesting question to ask when the algorithm will produce the first non-prime doodle.  The approach of choosing a 
region of maximal length, $p$, as the infinite region means that the doodle would contain a pair of finite regions with length no greater than $p$ that meet in a pair of distinct edges or vertices.  The smallest known example is the doodle shown in Figure 8.18, which has twenty crossings.

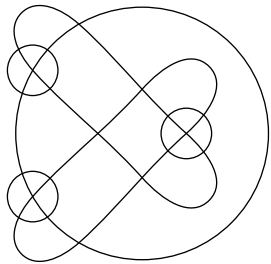
\diagram{A non-prime doodle produced by our algorithm}

The results of the search for reduced minimal prime doodles up to fourteen crossings are presented in the following table where $m$ is the number of components.
In each cell the first number gives the number of distinct doodles that are prime but not super-prime (3-connected but not 4-connected), the second number is the number of super-prime doodles.

The software used in the search may be found at www.layer8.co.uk/maths/doodles.

$$\vbox{ \offinterlineskip \halign{\strut \vrule \hfil \quad # 
                                \quad \hfil \vrule & \hfil \ # 
                                    \ \hfil \vrule & \hfil \ # 
                                    \ \hfil \vrule & \hfil \ # 
                                    \ \hfil \vrule & \hfil \ # 
                                    \ \hfil \vrule & \hfil \ # 
                                    \ \hfil \vrule \cr
\noalign{\hrule}
number of & \multispan4 prime minimal planar doodles \vrule \cr
\omit & \multispan4\hrulefill \cr
crossings & $m=1$ & $m=2$ & $m=3$ & $m=4$ \cr
\noalign {\hrule}
$6$ & $ $ & $ $ & $0,1$ & $ $ \cr
\noalign {\hrule}
$7$ & $ $ & $ $ & $ $ & $ $  \cr
\noalign {\hrule}
$8$ & $0,1$ & $ $ & $ $ & $ $ \cr
\noalign {\hrule}
$9$ & $1,0$ & $ $ & $ $  & $ $ \cr
\noalign {\hrule}
$10$ & $0,1$ & $0,1$ & $ $  & $ $ \cr
\noalign {\hrule}
$11$ & $1,0$ & $ $ & $1,1 $  & $ $ \cr
\noalign {\hrule}
$12$ & $1,2$ & $1,2$ & $1,1$  & $0,2$ \cr
\noalign {\hrule}
$13$ & $2,3$ & $4,1$ & $2,2$  & $2,0 $ \cr
\noalign {\hrule}
$14$ & $5,12$ & $8,14$ & $8,4$  & $0,1$ \cr
\noalign {\hrule}}}$$

\pagebreak
\section{Table: reduced prime and super prime doodles up to 14 crossings}

The following table lists all reduced minimal prime and super prime doodles up to 14 crossings.  The labelling takes the form $Xn^c_i$, where $X$ is $P$ or $S$ and indicates a prime or super prime doodle respectively, $n$ is the number of crossings, $c$ the number of components and $i$ an index.

Planar doodles with less than ten crossings

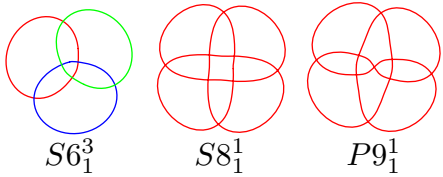
\picture{$S6^3_1$} \picture{$S8^1_1$} \picture{$P9^1_1$}

Planar doodles with ten crossings

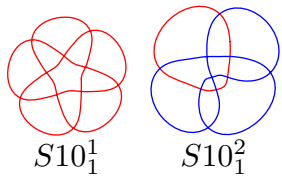
\picture{$S10^1_1$} \picture{$S10^2_1$} 

Planar doodles with eleven crossings

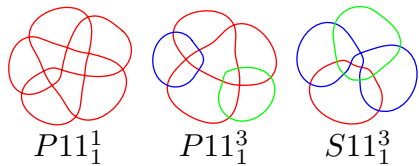
\picture{$P11^1_1$} \picture{$P11^3_1$} \picture{$S11^3_1$}

Planar doodles with twelve crossings

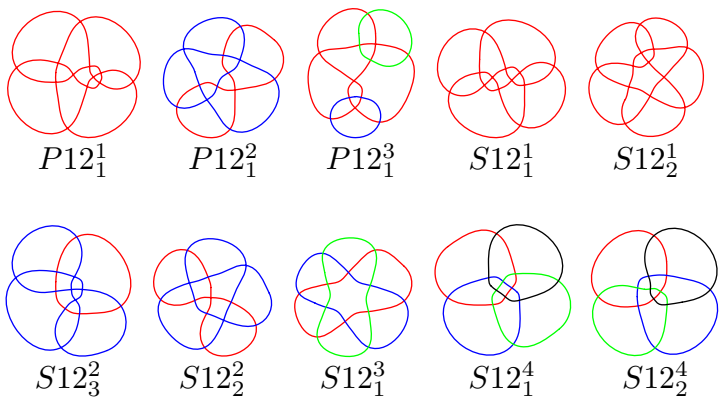
\picture{$P12^1_1$} \picture{$P12^2_1$} \picture{$P12^3_1$} \picture{$S12^1_1$} \picture{$S12^1_2$}

\picture{$S12^2_3$} \picture{$S12^2_2$} \picture{$S12^3_1$} \picture{$S12^4_1$} \picture{$S12^4_2$}

Planar doodles with thirteen crossings

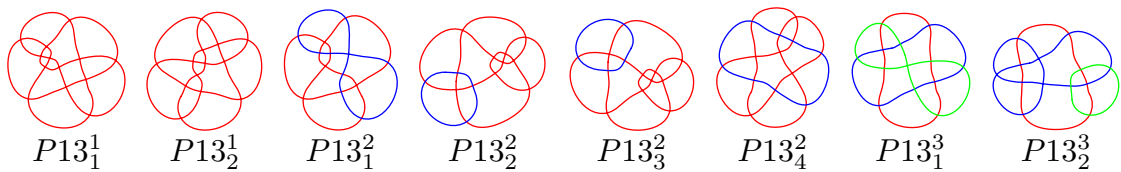
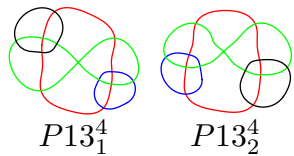
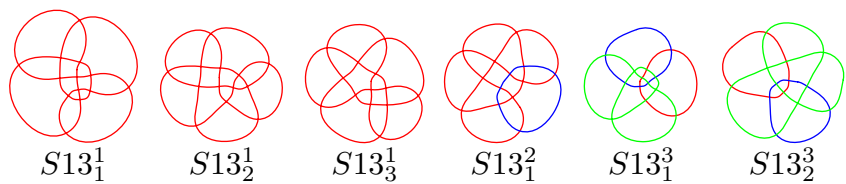
\picture{$P13^1_1$} \picture{$P13^1_2$} \picture{$P13^2_1$} \picture{$P13^2_2$} \picture{$P13^2_3$} \picture{$P13^2_4$} \picture{$P13^3_1$}  \picture{$P13^3_2$} 

\picture{$P13^4_1$}  \picture{$P13^4_2$} 

\picture{$S13^1_1$} \picture{$S13^1_2$} \picture{$S13^1_3$} \picture{$S13^2_1$} \picture{$S13^3_1$}  \picture{$S13^3_2$} 

Planar doodles with fourteen crossings

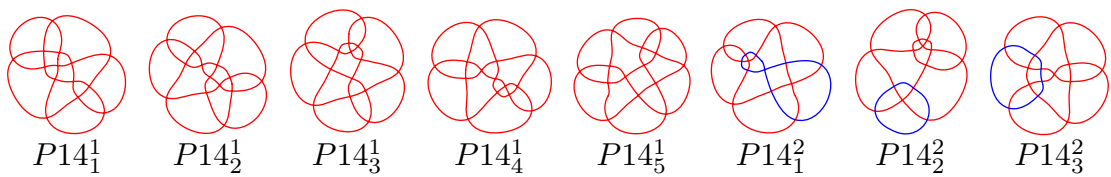
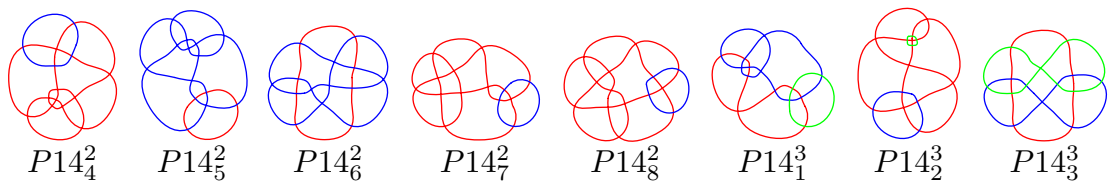
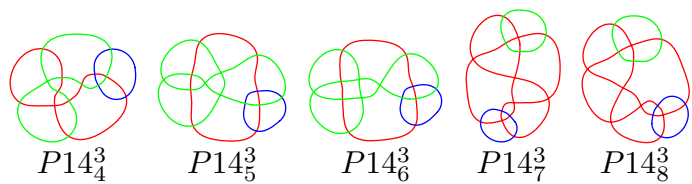
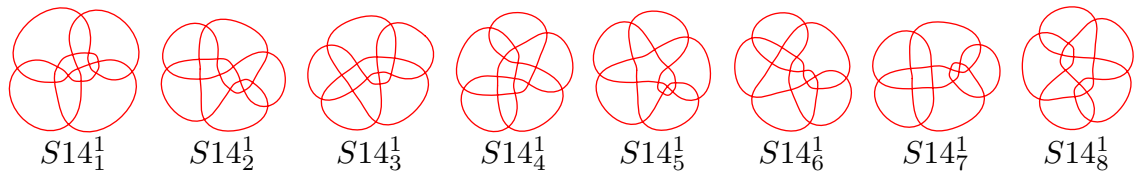
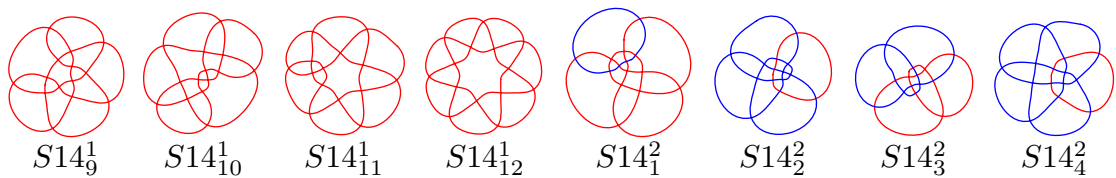
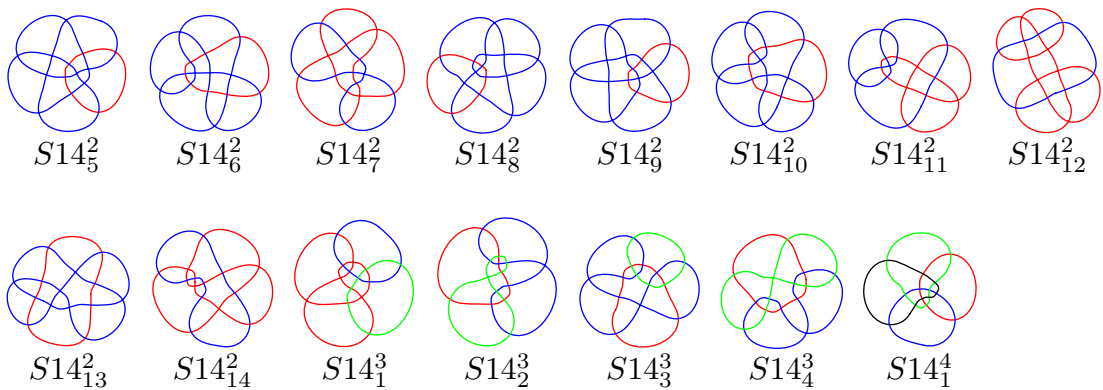
\picture{$P14^1_1$} \picture{$P14^1_2$} \picture{$P14^1_3$} \picture{$P14^1_4$} \picture{$P14^1_5$} \picture{$P14^2_1$} \picture{$P14^2_2$} \picture{$P14^2_3$} 

\picture{$P14^2_4$} \picture{$P14^2_5$} \picture{$P14^2_6$} \picture{$P14^2_7$} \picture{$P14^2_8$} \picture{$P14^3_1$} \picture{$P14^3_2$} \picture{$P14^3_3$} 

\picture{$P14^3_4$} \picture{$P14^3_5$} \picture{$P14^3_6$} \picture{$P14^3_7$} \picture{$P14^3_8$} 

\picture{$S14^1_1$} \picture{$S14^1_2$} \picture{$S14^1_3$} \picture{$S14^1_4$} \picture{$S14^1_5$} \picture{$S14^1_6$} \picture{$S14^1_7$} \picture{$S14^1_8$} 

\picture{$S14^1_9$} \picture{$S14^1_{10}$} \picture{$S14^1_{11}$} \picture{$S14^1_{12}$} \picture{$S14^2_1$} \picture{$S14^2_2$} \picture{$S14^2_3$} \picture{$S14^2_4$} 

\picture{$S14^2_5$} \picture{$S14^2_6$} \picture{$S14^2_7$} \picture{$S14^2_8$} \picture{$S14^2_9$} \picture{$S14^2_{10}$} \picture{$S14^2_{11}$} \picture{$S14^2_{12}$} 

\picture{$S14^2_{13}$} \picture{$S14^2_{14}$} \picture{$S14^3_1$} \picture{$S14^3_2$} \picture{$S14^3_3$} \picture{$S14^3_4$} \picture{$S14^4_1$}

\section{Last words}

We have described planar doodles and tabulated them up to 14 crossings. But as with all such investigations there are undiscovered pathways and pathways not completed. The application of hamiltonian circuits came rather late in our travails and we hope to answer the following questions in a later paper. \nl
1. Are there 2-connected or prime doodles with no hamiltonian circuit?\nl
2. Is there a simple proof of Tutte's result restricted to super prime doodles?\nl
3. Can the putative proof in 2. be applied to prime doodles?

\section{Bibliography}

[B] A. Bartholomew { Draw programme}. https://www.layer8.co.uk/maths/draw/

[BF] A. Bartholomew, R. Fenn. {\it Alexander and Markov Theorems for Generalized Knots, II Generalized Braids}, Journal of Knot Theory and its Ramifications Vol. 31, No. 8 (2022) 2240010.

[BFKK1] A. Bartholomew, R. Fenn, N. Kamada, S. Kamada. {\it Doodles on Surfaces}, arXiv:1612.08473v2, Journal of Knot Theory and its Ramifications Vol. 27 No 12 (2018), 1850071.

[BFKK2] A. Bartholomew, R. Fenn, N. Kamada, S. Kamada. {\it On Gauss codes of virtual doodles}, arXiv:1806.05885v1, special edition: Self-distributive system and quandle (co)homology theory in algebra and low-dimensional topology, Journal of Knot Theory and its Ramifications Vol. 27, No. 11 (2018) 1843013.

[F] R. Fenn, {\it Combinatorial Knot Theory}, to appear.

[FT] Roger Fenn and Paul Taylor, {\it Introducing doodles}, Topology of low-dimensional manifolds (Proc. Second Sussex Conf., Chelwood Gate, 1977), Lecture Notes in Math., vol. 722, Springer, Berlin, 1979, pp. 37–43.

 [FTur]
R.~Fenn, V.~Turaev, {\it Weyl Algebras and Knots}   { J. Geometry and
Physics} {\bf 57} (2007), pp.\ 1313--1324.

[HT] Hoste, J., and Thistlethwaite, M., { Knotscape 1.01}. https://web.math.utk.edu/~morwen/knotscape.html

[Kh] Mikhail Khovanov {\it Doodle groups} Trans. Amer. Math. Soc. 349 (1997), 2297-2315 

[T], W. Tutte.: {\it A theorem on planar graphs.} Trans. Amer. Math. Soc., 82(1956):
99-116.

\bye